\documentclass[11pt]{article}
\usepackage[a4paper, margin=0.75in, bmargin=0.85in]{geometry}
\usepackage[mathscr]{euscript}
\usepackage{amsmath,amsthm,mathrsfs}
\usepackage{amsfonts,amssymb}
\usepackage{graphicx}
\usepackage{breqn}
\usepackage{mathrsfs}
\usepackage{amsmath}
\usepackage{float}
\usepackage{caption}
\usepackage{chngcntr}
\usepackage{tikz-qtree}
\usetikzlibrary{arrows.meta}
\usepackage{subcaption}
\usepackage{epsfig,amsmath,amsthm,latexsym,amssymb,amsgen,graphicx}
\DeclareMathAlphabet{\mathpzc}{OT1}{pzc}{m}{it}
\usepackage{verbatim}
\usepackage{cite}
\usepackage{dsfont}
\usepackage{mathtools}
\DeclarePairedDelimiter{\ceil}{\lceil}{\rceil}
\newtheorem{thm}{Theorem}[section]
\newtheorem{cor}[thm]{Corollary}

\newtheorem{prop}[thm]{Proposition}

\newtheorem{defn}[thm]{Definition}
\newtheorem{rem}[thm]{Remark}
\newtheorem{exam}[thm]{Example}
\begin{document}
\title {\bf{Total Domination Index in Graphs}}
\author { \small Kavya. R. Nair \footnote{Corresponding author. Kavya. R. Nair, Department of Mathematics, National Institute of Technology, Calicut, Kerala, India. Email: kavyarnair@gmail.com} and M. S. Sunitha \footnote{sunitha@nitc.ac.in} \\ \small Department of Mathematics, National Institute of Technology, Calicut, Kerala, India-673601}
\date{}
\maketitle
\hrule
\begin{abstract}
 This paper introduces the concept of compliant vertices and compliant graphs, with a focus on the total domination degree (TDD) of a vertex in compliant graphs. The TDD is systematically calculated for various graph classes, including path graphs, cycles, book graphs, windmill graphs, wheel graphs, complete graphs, and complete bipartite graphs. The study explores inequalities involving TDD and defines total domination regular graphs. Furthermore, the TDD is analyzed in several graph operations such as union, join, composition, and corona, with a discussion on the property of the resulting graphs. The paper also examines the subdivision of complete graphs and degree splitting of path graphs. In the subsequent section, the total domination index (TDI) is introduced, and its values are calculated for different graph classes. The study concludes with bounds for the TDI across these graph classes.
\end{abstract}
\hrule
\section{Introduction}
Graphs are used to model pairwise relationships between objects. Graphs consist of vertices (also known as nodes) and edges that connect pairs of vertices.  Graph theory is an inevitable area of research because of its ability to depict and examine relationships, its universality, fundamental concepts, and widespread application in diverse domains. Graph theory has applications in the fields of computer science, operations research, social sciences, biology, and many more. Topological indices are numerical values that quantify a graph's structural traits or attributes. They are used in a variety of disciplines, including chemistry, network analysis, and algorithm creation. These indices offer insights, prediction power, and quantitative representations of graph topologies by capturing essential characteristics of graph connectivity or degrees. Among all the topological indices the most extensively studied and used index is the Wiener index \cite{W} which was developed by Harold Wiener in 1947 to compare the boiling point of some alkane isomers. Due to its widespread application numerous other topological indices were introduced which includes Randic Index (1975) \cite{r}, Zagreb Indices (1972) \cite{z1,z2}, Estrada Index (2000) \cite{e}, Gutman Index (1994) \cite{g}, Harary index (1993) \cite{h1,h2}, Balaban Index (1982) \cite{b1,b2} etc.  Domination in graphs is a fundamental concept that is vital due to its significance in wireless sensor networks, social networks, decision-making, and connection with other graph concepts. Ore and Berge's \cite{ore1962,berge1962} work on graph domination and related topics significantly contributed to the development and understanding of domination theory. Cockayne, Dawes, and Hedetniemi\cite{td} introduced the concept of total domination.
The notion of domination topological index in graph was introduced by Hanan et. al \cite{T5} in 2021. Motivated from the concept domination index in graphs was introduced by Kavya and Sunitha \cite{kav}. \\
This article is motivated by the idea of the domination index in graphs. The novel idea total domination index is introduced. Since it is not always possible to find a MTDS containing a particular vertex a new set of graphs called the compliant graphs are established.\\
The article begins with some of the basic definitions in graph theory. Section 3 introduces the concept of compliant vertices and compliant graphs. The idea of the total domination degree (TDD) of a vertex is defined in compliant graphs. The TDD is calculated for path graphs, cycles, book graphs, windmill graphs, wheel complete, and complete bipartite graphs. Some inequalities involving TDD are discussed. Total domination regular graphs are defined. The TDD is studied in operations of graphs such as union, join, composition, and corona. It is also established whether the resulting graphs are compliant or not. The subdivision of complete graphs and degree splitting of path graphs are also included in the study. Section 4 discusses the idea of TDI in graphs. The TDI of different graph classes is calculated using the results in Section 3. Bounds for TDI is also considered in the study. 
\section{Preliminaries}
The following fundamental definitions in graph theory are referred from \cite{W, T5, bal, char, k, thesis, tdtext, sub, degsp}
A graph $\mathcal{G}=(\mathcal{V},\mathcal{E})$, consists of a set $\mathcal{V}$ of vertices which is nonempty and a set $\mathcal{E}$ of edges which is a 2- element subset of $\mathcal{V}.$ An edge $e=ab$ means that $a$ and $b$ are adjacent or neighbors, and $a,b$ are end vertices of $e.$ Let $a$ be a vertex in $\mathcal{G}.$ The open neighborhood of $a$ is $\mathcal{N}(a)=\{b\,:\, \text{$ab$ is an edge}\}.$ The closed neighborhood of $a$ is $\mathcal{N}[a]=\{a\}\cup \mathcal{N}(a).$ The order $n$ and size $m$ of $\mathcal{G}$ are the number of vertices and edges respectively. The total number of edges incident at a vertex $a$ is called the degree of $a.$ The minimum and maximum degree of vertices in $\mathcal{G}$ is denoted as $\delta(\mathcal{G})$ and $\Delta(\mathcal{G})$ respectively. The distance between two vertices $a$ and $b$ in a connected graph $\mathcal{G}$ is the length of the shortest $a-b$ path. The eccentricity of a vertex $a$ is the maximum of distance between $a$ and any other vertex. The minimum and maximum eccentricity of $\mathcal{G}$ is the radius ($rad(\mathcal{G})$) and diameter ($diam(\mathcal{G})$) respectively. \\
Let $\mathcal{G}$ and $\mathcal{H}$ be two simple graphs. Then a graph isomorphism between $\mathcal{G}$ and $\mathcal{H}$ is a bijection $\theta:\mathcal{V}(\mathcal{G})\rightarrow \mathcal{V}(\mathcal{H})$ such that $a$ and $b$ are adjacent in $\mathcal{G}$ iff $\theta(a)$ and $\theta(b)$ are adjacent in $\mathcal{H}.$ A graph $\mathcal{H}$ is a subgraph of $\mathcal{G}$ if $\mathcal{V}(\mathcal{H})\subseteq \mathcal{V}(\mathcal{G})$ and $\mathcal{E}(\mathcal{H})\subseteq \mathcal{E}(\mathcal{H}).$ If $\mathcal{V}(\mathcal{G})=\mathcal{V}(\mathcal{H}),$ then $\mathcal{H}$ is a spanning subgraph. Now, $\mathcal{T}$ is a spanning tree of $\mathcal{G},$ is $\mathcal{T}$ is both a spanning subgraph and a tree. 
Let $\mathcal{G}_1=(\mathcal{V}_1,\mathcal{E}_1)$ and $\mathcal{G}_2=(\mathcal{V}_2,\mathcal{E}_2)$ be two graphs. Then the union $\mathcal{G}=\mathcal{G}_1\cup \mathcal{G}_2$ is a graph with $\mathcal{V}=\mathcal{V}_1\cup \mathcal{V}_2$ and $\mathcal{E}=\mathcal{E}_1\cup \mathcal{E}_2.$ The join $\mathcal{G}=\mathcal{G}_1+\mathcal{G}_2$ is a graph with $\mathcal{V}=\mathcal{V}_1\cup \mathcal{V}_2$ and $\mathcal{E}=\mathcal{E}_1\cup \mathcal{E}_2\cup \mathcal{E}',$ where $\mathcal{E}'$ is the set of edges that join each vertex of $\mathcal{G}_1$ to every vertex of $\mathcal{G}_2.$ Let $\mathcal{G}= \mathcal{G}_1\circ \mathcal{G}_2$ be the composition of $\mathcal{G}_1$ and $\mathcal{G}_2.$ Then $\mathcal{V}(\mathcal{G})=\mathcal{V}_1\times \mathcal{V}_2.$ Two vertices $(a_1,b_1)$ and $(a_2,b_2)$ are adjacent in composition $\mathcal{G}= \mathcal{G}_1\circ \mathcal{G}_2,$ iff  $a_1$ is adjacent to $a_2$
in $\mathcal{G}_1$ or $a_1=a_2$ and $b_1$ is adjacent to $b_2$ in $\mathcal{G}_2.$ The corona $\mathcal{G}=\mathcal{G}_1\odot\mathcal{G}_2$ is obtained by taking $\lvert\mathcal{V}_1\rvert$ copies of $\mathcal{G}_2$ and joining $i^{th}$ vertex if $\mathcal{G}_1$ to every vertex in the $i^{th}$ copy of $\mathcal{G}_2.$  \\
A wheel graph $\mathcal{W}_n$ is the join of $\mathcal{K}_1$ and cycle $\mathcal{C}_n.$ A Windmill graph $\mathcal{W}d(r,s)$ is obtained by taking $s$ copies of $\mathcal{K}_r$ with a common vertex. The common vertex is called center vertex of $\mathcal{W}d(r,s)$. The Cartesian product of a star graph $\mathcal{S}_{n+1}$ and path $\mathcal{P}_2$ is a book graph $\mathcal{B}_n.$ If the vertex in $\mathcal{S}_{n+1}$ with degree $n$ is $a$ and vertices of $\mathcal{P}_2$ are $b_1$ and $b_2$, then the centre vertex of $\mathcal{B}_n$ are $(a,b_1)$ and $(a,b_2).$
Now, consider $\mathcal{P}_3.$ For $s=2,3,...$ identify the roots of $s$ copies of $\mathcal{P}_3$ to obtain a rooted tree $\mathscr{B}_s$, and the root of $\mathscr{B}$ is the vertex obtained by identifying copies of $\mathcal{P}_3.$
Let $t\geq 2,$ and $\beta_1,\beta_2,...\beta_t\in \{\mathscr{B}_1,\mathscr{B}_2,...\}.$ A Kragujevac tree, $\mathcal{T}$, is a tree that includes a vertex say $a$ with a degree of $t$. The vertex $a$ is adjacent to the roots of $\beta_1,\beta_2,...\beta_t.$ Vertex $a$ is referred to as the central vertex of $\mathcal{T}$. The subgraphs $\beta_1,\beta_2,...\beta_t$ are called the branches of $\mathcal{T}$. The order of $\mathscr{B}_s$ is $2s+1.$ Hence a Kragujevac tree with $t$ branch defined as $\beta_i \cong \mathscr{B}_{s_i}$, $i=1,2,...,t$ has order $n=1 +\sum\limits_{i=1}^t (2s_i + 1).$ 
The subdivision S($\mathcal{G}$) acts on $\mathcal{G}$ by replacing each of its edges by a path of length two. Now, let $V(\mathcal{G})=A_1\cup A_2\cup...,\cup A_k\cup S ,$ where each $A_i$ is collection vertices having same degree and $S$ is the set of remaining vertices. Then degree splitting of $\mathcal{G}$ is obtained by adding vertices $a_1,a_2,...,a_k$ and joining $a_i$ to each vertex in $A_i, i=1,2,...,k.$ \\
Let $\mathcal{D}\subseteq \mathcal{V},$ then $\mathcal{D}$ is a dominating set (DS) if for every vertex  $b\in \mathcal{V}\setminus \mathcal{D}$ there exists $a\in \mathcal{D}$ such that $a$ is adjacent to $b,$ i.e, every vertex in $\mathcal{G}$ is either in $\mathcal{D}$ or is a neighbor of some vertex in $\mathcal{D}.$ The minimum cardinality of a DS in $\mathcal{G}$ is called domination number of $\mathcal{G}$ denoted by $\gamma(\mathcal{G})$ or simply $\gamma.$ The DS with the least number of vertices is called a minimum DS. A DS is a minimal dominating set (MDS) if it does not properly contain any other DS. 
Let $\mathcal{G}$ be a graph without isolated vertices. A total dominating set (TDS), $\mathcal{S}$ of a graph $\mathcal{G}$ is a set of vertices of $\mathcal{G}$ such that every vertex is adjacent to a
 vertex in $\mathcal{S}$, i.e, $\mathcal{N}(\mathcal{S})=V$. If no proper subset of
 $\mathcal{S}$ is a TDS , then $\mathcal{S}$ is a minimal TDS (MTDS) of $\mathcal{G}$. The minimum cardinality of TDS of $\mathcal{G}$ is the total domination number (TDN), $\gamma_t(\mathcal{G})$ (or simply $\gamma_t$) of $\mathcal{G}.$ The maximum cardinality of MTDS in $\mathcal{G}$ is called the upper total domination number $\Gamma_t(\mathcal{G}).$ A vertex is a good vertex if it is contained in any minimum TDS. A graph is $\gamma_t-$excellent if every vertex of the graph is a good vertex. 
\section{Total domination degree in graphs}
The novel notion of total domination degree, compliant vertex and compliant graphs are introduced in the section. Through an illustration, it is established that it is not always possible to find a MTDS containing a particular vertex. Thus a new graph class is introduced namely compliant graphs.
\begin{figure}
    \centering
    \includegraphics[width=0.5\linewidth]{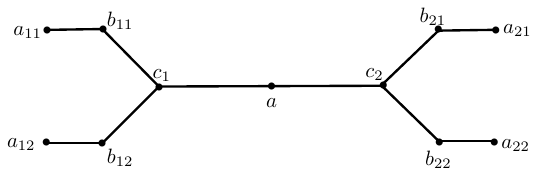}
    \caption{Illustration for TDD of a vertex.}
    \label{f1}
\end{figure}
\begin{exam}\label{e1}
Consider the graph in Figure \ref{f1}. Take the vertex $c_1.$ The set $\{c_1,b_{12},b_{11},c_2,b_{21},b_{22}\}$ is a MTDS containing $c_1.$ The same set can be considered as a MTDS containing $b_{12},b_{11},c_2,b_{21}$ and $b_{22}.$ Now take the end vertices. The set $\{a_{11},b_{11},a_{12},b_{12},c_2,b_{21},b_{22}\}$ is a MTDS containing $a_{11}.$ The same set can be taken as the MTDS containing $a_{12}$ also. Similarly for $a_{21}$ and $a_{22}$ the set $\{a_{21},b_{21},a_{22},b_{22},c_1,b_{11},b_{12}\}$ can be considered as the required MTDS. For the vertex $a,$ either $c_1$ or $c_2$ is required for the total domination. Suppose $c_1$ is included in the set, then $b_{11}$ and $b_{12}$ are required to totally dominate $a_{11}$ and $a_{12}$ respectively. Now $a_{21},b_{21},a_{22},b_{22}$ can be included to totally dominate the remaining vertices. So, set $\{a,a_{21},b_{21},a_{22},b_{22},c_1,b_{11},b_{12}\}$ is a TDS containing $a.$ But this set already contains a TDS, hence, it is not minimal. Therefore it can be concluded that there does not exist a MTDS containing $a$ in the graph.     
\end{exam}
From \ref{e1}, the definition for compliant vertex, non-compliant vertex, and compliant graphs can be derived as follows:
\begin{defn}
    Let $\mathds{G}$ be a graph and let $a$ be a vertex in $\mathds{G}$. Then $a$ is called a compliant vertex if there exists a MTDS containing $a$ in $\mathds{G}$. Otherwise, $a$ is called a non- compliant vertex. 
\end{defn}
\begin{defn}
    A graph $\mathds{G}$ is called compliant graph if every vertex in $\mathds{G}$ is a compliant vertex.  
\end{defn}
From Example \ref{e1} it can be concluded that Kragujevac trees are non compliant graphs with one non- compliant vertex.
\begin{figure}
    \centering
    \includegraphics[width=0.5\linewidth]{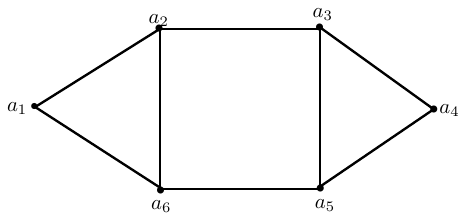}
    \caption{A compliant graph}
    \label{f2}
\end{figure}
\begin{exam}\label{e2}
For the graph in Figure \ref{f2}, every vertex is a compliant vertex. Here sets $\{a_1,a_2,a_5,a_4\}$ and $\{a_1,a_6,a_3,a_4\}$ can be taken as required sets for all vertices. Also it can be seen that $\{a_2,a_3\}$ is also a MTDS containing $a_2$ and $a_3$ of cardinality 2. Hence for $a_2$ and $a_3$ the MTDS containing each of them with least cardinality is $\{a_2,a_3\}.$
\end{exam}
Now, from Example \ref{e2} definition for TDD of a vertex in a compliant graph is derived.
\begin{defn}
Let $\mathds{G}$ be a compliant graph and $a$ be a vertex in $\mathds{G}$. The total domination degree (TDD) of vertex $a$ is the minimum number of vertices in a MTDS containing $a,$ i.e,
    $$\mathpzc{d}_{td}(a)=min \{|S|: \text{ $S$ is a MTDS containing $a$}\}$$
    The minimum and maximum TDD of $\mathds{G},$ denoted as $\delta_{td}(\mathds{G})$ and $\Delta_{td}(\mathds{G})$ respectively, are defined as,
    $$_{td}\delta(\mathds{G})= min\,\{\mathpzc{d}_{td}(a):\, a\in \mathbb{V}(\mathds{G})\}$$
    $$_{td}\Delta(\mathds{G})= max\,\{\mathpzc{d}_{td}(a):\, a\in \mathbb{V}(\mathds{G})\}.$$
\end{defn}
\begin{exam}
    For the compliant graph in Figure \ref{f2}, the TDD of vertices $a_1$ and $a_4$ is 4. And for all other vertices, the TDD is 2. Hence $\delta_{td}(\mathds{G})=2$ and $\Delta_{td}(\mathds{G})=4.$
\end{exam}
Now, some classes of graphs are classified as compliant graphs.
\subsection{Total domination index of Paths $\mathds{P}_n$ $n\neq 4,7$}
Consider $\mathds{P}_4$ and label the vertices as $a_1,a_2,a_3$ and $a_4.$ For $a_2$ and $a_3$ the set $\{a_2,a_3\}$ can be taken as the MTDS containing each of them. Now, take one of the end vertices say $a_1.$ To totally dominate $a_1,$ $a_2$ is needed. Include $a_3$ to totally dominate every vertex. But the set $\{a_1,a_2,a_3\}$ is not a MTDS. Hence there does not exist a MTDS containing $a_1.$ Similarly there does not exist a MTDS containing $a_4.$ Hence in the case of $\mathds{P}_4$ the end vertices are non- compliant vertices. And the support vertices are compliant vertices. \\
Similarly consider path $\mathds{P}_7$ and label the vertices as $a_1,a_2,...,a_7.$ The sets $\{a_1,a_2,a_5,a_6\}$ and $\{a_2,a_3,a_6,a_7\}$ can be taken as MTDS containing each of $a_1,a_2,a_3,a_5,a_6$ and $a_7.$ Hence $a_1,a_2,a_3,a_5,a_6$ and $a_7$ are compliant vertices. Now, take vertex $a_4.$ For total domination either $a_3$ or $a_5$ is required. Suppose $a_5$ is included. Then $a_6$ is also needed to totally dominate $a_7.$ Also, $a_1$ and $a_2$ must be included to totally dominated the entire vertex set. Hence the obtained set is $\{a_1,a_2,a_4,a_5,a_6\}$ which is not minimal. Likewise, it can be concluded there does not exist a MTDS containing $a_4$ in $\mathds{P}_7.$ Hence $\mathds{P}_7$ is not compliant graph and $a_4$ is the non- compliant vertex.\\\\
For a path $\mathds{P}_n,$ label the vertices as $a_1,a_2,...a_n.$
\begin{thm}\cite{thesis}\label{t1}
 Path $\mathds{P}_n$ is $\gamma_t-$ excellent for $n \equiv 2$ (mod $4$).   
\end{thm}
\begin{thm}\cite{tdtext}\label{t2}
    For $n \geq 3, $ 
    $$\gamma_t(\mathds{P}_n)=\gamma_t(\mathds{C}_n)=\begin{cases}
      \frac{n}{2} \quad &\text{ if  } n\equiv 0 (mod4)\\
      \frac{n+1}{2} \quad &\text{ if  } n\equiv 1,3 (mod 4)\\
      \frac{n}{2}+1 \quad &\text{ if  } n\equiv 2 (mod  4)
    \end{cases}$$
\end{thm}
\begin{thm}\label{t3}
    The path graph $\mathds{P}_n$ is a compliant graph for $n \neq 4,7.$
\end{thm}
\begin{proof}
    The case of $\mathds{P}_2$ and $\mathds{P}_3$ can be proved easily. For $\mathds{P}_5,$ set $\{a_2,a_3,a_4\}$ and $\{a_1,a_2,a_4,a_5\}$ can be taken as required MTDS containing each of the vertices. Similarly for $\mathds{P}_6,$ the sets $\{a_1,a_2,a_5,a_6\},$ $\{a_2,a_3,a_5,a_6\}$ and $\{a_1,a_2,a_4,a_5\}$ can be considered as required MTDS containing each of the vertices. Hence $\mathds{P}_5$ and $\mathds{P}_6$ are compliant vertices. Now, for $n\geq 8,$ classify $\mathds{P}_n$ to 4 cases. \\
    \textbf{Case 1}: $n\equiv 2 (\textrm{ mod } 4).$\\
    Let $n=4k+2, k>1, k\in \mathbb{N}.$ From Theorem \ref{t1}, it is clear that every vertex in $\mathds{P}_n$ belongs to a minimum TDS which is a MTDS. Hence, clearly $\mathds{P}_n$ is compliant for $n\equiv 2 (\textrm{ mod } 4).$ Since every vertex belongs to a minimum TDS, by Theorem \ref{t2} the cardinality of every minimum TDS is $\gamma_t=2k+2.$ \\
    For vertices $a_i,$ $i\equiv 1,2 (\textrm{ mod }4)$ take the set $\{\{a_1,a_2\},\{a_5,a_6\},\{a_9,a_{10}\},...,\{a_{4k-3}a_{4k-2}\},\{a_{4k+1},a_{4k+2}\}\}$ which contains 2k+2 vertices and totally dominates the entire vertex set. Hence the given set is minimum TDS and therefore is MTDS containing each of the vertices $a_i,$ $i\equiv 1,2 (\textrm{ mod }4)$.\\
    For vertices $a_i,$ $i\equiv 0 (\textrm{ mod }4)$  consider the set $\{\{a_1,a_2\},\{a_4,a_5\},\{a_8,a_9\},...,\{a_{4k-4}a_{4k-3}\}, \{a_{4k},a_{4k+1}\}\}$ is a TDS containing 2k+2 vertices. Hence the set is minimum TDS, therefore a MTDS containing each of the vertices $a_i,$ $i\equiv 0 (\textrm{ mod }4)$.
    Similarly for vertices $a_i,$ $i\equiv 3 (\textrm{ mod }4),$ consider the set $\{\{a_2,a_3\},\{a_6,a_7\},\{a_{10},a_{11}\},$ $...,\{a_{4k-2}a_{4k-1}\}, \{a_{4k+1},a_{4k+2}\}\}$ is a TDS containing 2k+2 vertices. Hence the set is minimum TDS, therefore a MTDS containing each of the vertices $a_i,$ $i\equiv 3 (\textrm{ mod }4)$.\\
    \textbf{Case 2:} $n\equiv 3 (\textrm{ mod } 4).$\\
    Let $n=4k+3, k>1, k\in \mathbb{N}.$
    For $n\equiv 3 (\textrm{ mod } 4),$ by Theorem \ref{t2} $\gamma_t=2k+2.$ Hence any TDS of cardinality $2k+2$ is a minimum TDS and therefore a MTDS.\\
    For vertices $a_i,$ $i\equiv 2,3 (\textrm{ mod }4),$ consider the set $\{\{a_2,a_3\},\{a_6,a_7\},\{a_{10},a_{11}\},...,\{a_{4k-2}a_{4k-1}\},$\\$ \{a_{4k+2},a_{4k+3}\}\}$ is a TDS containing 2k+2 vertices. Hence the set is minimum TDS, therefore a MTDS containing each of the vertices $a_i,$ $i\equiv 2,3 (\textrm{ mod }4)$.\\
    For vertices $a_i,$ $i\equiv 1 (\textrm{ mod }4),$ consider the set $\{\{a_1,a_2\},\{a_5,a_6\},\{a_{9},a_{10}\},$ $...,\{a_{4k-3}a_{4k-2}\}, \{a_{4k+1},a_{4k+2}\}\}$ is a TDS containing 2k+2 vertices. Hence the set is minimum TDS, therefore a MTDS containing each of the vertices $a_i,$ $i\equiv 1 (\textrm{ mod }4)$.\\
    For vertices $a_i,$ $i\equiv 0 (\textrm{ mod }4),$ consider the set $\{\{a_2,a_3,a_4\},\{a_7,a_{8}\},\{a_{11},a_{12}\}...,\{a_{4k-1}a_{4k}\}, \{a_{4k+2},a_{4k+3}\}\}$ is a TDS containing 2k+3 vertices. If any vertex in an element of the set is deleted the other vertex in the same element is not totally dominated. Similarly, if $a_2$ or $a_4$ are deleted $a_1$ and $a_5$ are not totally dominated respectively. Likewise deleting $a_3$ disturbs the total domination of $a_2$ and $a_4.$ Hence the given set is a MTDS containing each of the vertices $a_i,$ $i\equiv 0 (\textrm{ mod }4)$.\\
    \textbf{Case 3:} $n\equiv 0 (\textrm{ mod } 4).$\\
    Let $n=4k, k>1, k\in \mathbb{N}.$ For $n\equiv 0 (\textrm{ mod } 4),$ by Theorem \ref{t2} $\gamma_t=2k.$ \\
    For vertices $a_i,$ $i\equiv 2,3 (\textrm{ mod }4),$ consider the set $\{\{a_2,a_3\},\{a_6,a_7\},\{a_{10},a_{11}\},$ $...,\{a_{4k-2},a_{4k-1}\}\}$ is a TDS containing 2k vertices. Hence the set is minimum TDS, therefore a MTDS containing each of the vertices $a_i,$ $i\equiv 2,3 (\textrm{ mod }4)$.\\
    For vertices $a_i,$ $i\equiv 0 (\textrm{ mod }4),$
    the set $\{\{a_2,a_3,a_4\},\{a_7,a_8\},\{a_{11},a_{12}\},...\{a_{4k-5},a_{4k-4}\},\{a_{4k-1},a_{4k}\}\}$ is a TDS containing 2k+1 vertices. If any $a_2$ or $a_4$ is deleted, then $a_1$ or $a_5$ is not totally dominated. And deleting $a_3$ disturbs the total domination of $a_2$ and $a_4$. Also if any vertex in any other element of the set is deleted the other vertex in the same element is not totally dominated. Hence the set is MTDS containing $a_i,$ $i\equiv 0 (\textrm{ mod }4)$.\\
    For vertices $a_i,$ $i\equiv 1 (\textrm{ mod }4),$
    the set $\{\{a_1,a_2\},\{a_5,a_6\},\{a_{9},a_{10}\},...\{a_{4k-7},a_{4k-6}\},\{a_{4k-3},a_{4k-2},a_{4k-1}\}\}$ is a TDS containing 2k+1 vertices. If any $a_{4k-3}$ or $a_{4k-1}$ is deleted, then $a_{4k-4}$ or $a_{4k}$ is not totally dominated. And deleting $a_{4k-2}$ disturbs the total domination of $a_{4k-3}$ and $a_{4k-1}$. Also if any vertex in any other element of the set is deleted the other vertex in the same element is not totally dominated. Hence the set is MTDS containing $a_i,$ $i\equiv 1 (\textrm{ mod }4)$.\\
    \textbf{Case 4:} $n\equiv 1 (\textrm{ mod } 4).$\\
    Let $n=4k+1, k>1, k\in \mathbb{N}.$ For $n\equiv 0 (\textrm{ mod } 4),$ by Theorem \ref{t2} $\gamma_t=2k+1.$ \\
    For vertices $a_i,$ $i\equiv 2,3 (\textrm{ mod }4),$ consider the set $\{\{a_2,a_3\},\{a_6,a_7\},\{a_{10},a_{11}\},$ $...,\{a_{4k-6},a_{4k-5}\},$ \\$\{a_{4k-2},a_{4k-1},a_{4k}\}\}$ is a TDS containing 2k+1 vertices. Hence the set is minimum TDS, therefore a MTDS containing each of the vertices $a_i,$ $i\equiv 2,3 (\textrm{ mod }4)$.\\
    For vertices $a_i,$ $i\equiv 0 (\textrm{ mod }4),$ consider the set $\{\{a_2,a_3,a_4\},\{a_7,a_8\},\{a_{11},a_{12}\},$ $...,\{a_{4k-5},a_{4k-4}\},\{a_{4k-1},a_{4k}\}\}$ is a TDS containing 2k+1 vertices. Hence the set is minimum TDS, therefore a MTDS containing each of the vertices $a_i,$ $i\equiv 0 (\textrm{ mod }4)$.\\ 
    For vertices $a_i,$ $i\equiv 1 (\textrm{ mod }4),$
    the set $\{\{a_1,a_2\},\{a_5,a_6\},\{a_{9},a_{10}\},...\{a_{4k-3},a_{4k-2}\},\{a_{4k-1},a_{4k}\}\}$ is a TDS containing 2k+2 vertices. If any vertex in an element of the set is deleted the other vertex in the same element is not totally dominated. Hence the set is MTDS containing $a_i,$ $i\equiv 1(\textrm{ mod }4)$.\\
\end{proof}
\begin{thm}\label{t4}
    For $\mathds{P}_n,$ $n\neq 4,7$ 
    $$\mathpzc{d}_{td}(a_i)=\begin{cases}
    \frac{n}{2}+1 &\text{ if }
n \equiv 2 ( mod 4)\\
\frac{n+1}{2} & \text{ if } n\equiv 1(mod 4) \text{ and } i\not\equiv 1(mod 4)\\
\frac{n+1}{2}+1 & \text{ if } n\equiv 1(mod 4) \text{ and } i\equiv 1(mod 4)\\
\frac{n+1}{2} & \text{ if } n\equiv 3(mod 4) \text{ and } i\not\equiv 0(mod 4)\\
\frac{n+1}{2}+1 & \text{ if } n\equiv 3(mod 4) \text{ and } i\equiv 0(mod 4)\\
\frac{n}{2} & \text{ if } n\equiv 0(mod 4) \text{ and } i\equiv 2,3(mod 4)\\
\frac{n}{2}+1 & \text{ if } n\equiv 0(mod 4) \text{ and } i\equiv 0,1(mod 4)
    \end{cases}$$
\end{thm}
\begin{proof}
From Theorem \ref{t3} it is clear that $\mathbb{P}_n,$ $n\neq 4,7$ is a compliant graph. Now, it is required to find the minimum number of vertices in the MTDS containing each vertex. \\
\textbf{Case 1:} $n\equiv 2(mod 4).$\\
From Theorem \ref{t2} the TDN $\gamma_t$ for this case is $\frac{n}{2}+1.$ Hence TDS containing $\frac{n}{2}+1$ elements is minimum and hence minimal. And in proof of Theorem \ref{t3} it is proved that each vertex belongs to a MTDS containing $2k+2$ vertices where $n=4k+2.$ Hence $\mathpzc{d}_{td}(a_i)=\frac{n}{2}+1$ for all $i=1,2,...,n.$\\
\textbf{Case 2:} $n\equiv 1(mod 4).$\\
The TDN $\gamma_t$ for this case is $\frac{n+1}{2}.$ From the proof of theorem \ref{t3} it is clear that vertices $a_i$ such that $i\not\equiv1(mod4)$ is contained in MTDS containing $2k+1$ vertices where $n=4k+1.$ Therefore, $\mathpzc{d}_{td}(a_i)=\frac{n+1}{2}$ for $i\not\equiv1(mod 4).$ Now for vertices $a_i$ such that $i\equiv 1(mod 4)$ it is already proved that there exists MTDS containing $a_i$ with $2k+2$ vertices where $n=4k+1.$ Hence it is enough to prove that there does not exist a MTDS containing $a_i$ with $2k+1$ vertices. For that, consider the path $\mathbb{P}_n,n=4k+1$ and split the path into two paths: the first path is $a_1a_2a_3...a_{i-2}$ and the second path is $a_{i-1}a_i...a_{4k+1}.$ It is required to find a MTDS containing $a_i,i=4m+1,m<k.$ Hence the second path is taken as a path in which $a_i$ is a support vertex. The first path is a path of length $4m-1$ and the second path is of length $4(k-m)+2.$ From Theorem \ref{t2} at least $2m$ vertices are required to totally dominate the first path and a TDS containing $a_i$ with at least $2(k-m)+2$ vertices are required to totally dominate the second path. Hence in total, a TDS containing $a_i$ has at least $2k+2=\frac{n+1}{2}+1$ vertices in $\mathbb{P}_n.$ Therefore, $\mathpzc{d}_{td}(a_i)=\frac{n+1}{2}+1$ for $i\equiv1(mod 4).$ \\
\textbf{Case 3:} $n\equiv 0(mod 4).$\\
The TDN $\gamma_t$ for this case is $\frac{n}{2}.$ From the proof of Theorem \ref{t3} it is clear that vertices $a_i$ such that $i\not\equiv0,1(mod4)$ is contained in MTDS containing $2k$ vertices where $n=4k.$ Therefore, $\mathpzc{d}_{td}(a_i)=\frac{n}{2}$ for $i\not\equiv0,1(mod 4).$ Now for vertices $a_i$ such that $i\equiv 0(mod 4)$ it is already proved that there exists MTDS containing $a_i$ with $2k+1$ vertices where $n=4k.$ Hence it is enough to prove that there does not exist a MTDS containing $a_i$ with $2k$ vertices. For that consider the path $\mathbb{P}_n,n=4k$ and split the path into two paths: the first path is $a_1a_2a_3...a_{i-2}a_ia_{i+1}$ and the second path is $a_{i+2}a_{i+3}...a_{4k}, i\neq 4k.$ It is required to find a MTDS containing $a_i,i=4m,m<k.$ Hence the first path is taken as a path in which $a_i$ is a support vertex. The first path is a path of length $4m+1$ and the second path is of length $4(k-m)+1.$ From Theorem \ref{t2} at least $2(k-m)$ vertices are required to totally dominate the second path and a TDS containing $a_i$ with at least $2m+1$ vertices are required to totally dominate the first path. a TDS containing $a_i$ has at least $2k+1=\frac{n}{2}+1$ vertices in $\mathbb{P}_n.$ For $a_i,i=4k,$ it is clear that 2 vertices namely $a_{4k}$ and $a_{4k-1}$ are required to totally dominate $a_{4k-2}a_{4k-1}a_{4k}.$ The remaining vertices can be considered as a path $a_1a_2...a_{4k-3}$ of length $4k-3.$ And by Theorem \ref{t2} at least $2k-1$ vertices are required to totally dominate the path. Hence a TDS containing $a_i$ has at least $2k+1=\frac{n}{2}+1$ vertices in $\mathbb{P}_n.$ Therefore, $\mathpzc{d}_{td}(a_i)=\frac{n}{2}+1$ for $i\equiv0(mod 4).$
By the same argument for $a_i,i=4m+1,m<k,$ path $\mathbb{P}_n$ is split into two paths: first path is $a_1a_2...a_{4m-1}$ and the second path is $a_{4m}a_{4m-1}...a_{4k}.$ Hence in total at least $2k+1$ vertices are required to totally dominate $\mathbb{P}_n.$ Therefore, $\mathpzc{d}_{td}(a_i)=\frac{n}{2}+1$ for $i\equiv1(mod 4).$\\
\textbf{Case 4:} $n\equiv 3(mod 4).$\\
The TDN $\gamma_t$ for this case is $\frac{n+1}{2}.$ From the proof of Theorem \ref{t3} it is clear that vertices $a_i$ such that $i\not\equiv0(mod4)$ is contained in MTDS containing $2k+2$ vertices where $n=4k+3.$ Therefore, $\mathpzc{d}_{td}(a_i)=\frac{n+1}{2}$ for $i\not\equiv0(mod 4).$ For vertices $a_i, i\equiv 0(mod 4),$ split the path into two: first path is $a_1a_2...a_{4m-2}$ of length $4m-2$ and the second path is $a_{4m-1}a_{4m}...a_{4k+3}$ of length $4((k+1)-m)+1$ where $i=4m.$ Hence, at least $2m$ and $2(k-m)+3$ vertices are required to totally dominate the first and second paths respectively. Therefore TDS containing $a_i$ consists of at least $2k+3$ vertices. Hence $\mathpzc{d}_{td}(a_i)=\frac{n+1}{2}+1$ for $i\equiv0(mod 4).$ 
\end{proof}
\begin{thm}\label{t5}
   Cycle $\mathbb{C}_n,n\geq 3$ is a compliant graph and $$\mathpzc{d}_{td}(a_i)=\begin{cases}
     \frac{n}{2} \quad &\text{ if  } n\equiv 0 (mod4)\\
      \frac{n+1}{2} \quad &\text{ if  } n\equiv 1,3 (mod 4)\\
      \frac{n}{2}+1 \quad &\text{ if  } n\equiv 2 (mod  4)   
   \end{cases}$$ 
\end{thm}
\begin{proof}
    From Theorem \ref{t2} the TDN of a cycle can be determined. Hence it is enough to prove that every vertex in $\mathbb{C}_n$ is contained in a TDS containing $\gamma_t$ vertices. Consider 4 cases:\\
    \textbf{Case 1:} $n\equiv0(mod 4).$\\
    Let $n=4k, k\in \mathbb{N}.$ For vertices of the form $a_i,i\equiv 1,2(mod 4)$ take the set $\{\{a_1,a_2\},\{a_5,a_6\},$ $...,\{a_{4k-3},a_{4k-2}\}\}$ which is a TDS containing $2k=\frac{n}{2}$ vertices. Similarly for vertices $a_i,i\equiv 0,3(mod 4)$ consider the set $\{\{a_3,a_4\},\{a_7,a_8\},...\{a_{4k-1},a_{4k}\}\}$ which is a TDS containing $2k=\frac{n}{2}$ vertices. Hence $\mathpzc{d}_{td}(a_i)=\frac{n}{2}.$\\
    \textbf{Case 2:} $n\equiv1(mod 4).$\\
    Let $n=4k+1, k\in \mathbb{N}.$ For vertices of the form $a_i,i\equiv 1,2(mod 4)$ take the set $\{\{a_1,a_2\},\{a_5,a_6\},$\\$...,\{a_{4k-3},a_{4k-2}\},\{a_{4k-1}\}\}$ which is a TDS containing $2k+1=\frac{n+1}{2}$ vertices. Similarly for vertices $a_i,i\equiv 0,3(mod 4)$ consider the set $\{\{a_3,a_4\},\{a_7,a_8\},$\\$...,\{a_{4k-1},a_{4k}\}.\{a_{4k+1}\}\}$ which is a TDS containing $2k+1=\frac{n+1}{2}$ vertices. Hence $\mathpzc{d}_{td}(a_i)=\frac{n+1}{2}.$\\
    \textbf{Case 3:} $n\equiv3(mod 4).$\\
    Let $n=4k+3, k\in \mathbb{N}\cup\{0\}.$ For vertices of the form $a_i,i\equiv 1,2(mod 4)$ take the set $\{\{a_1,a_2\},\{a_5,a_6\},$\\$...,\{a_{4k+1},a_{4k+2}\}\}$ which is a TDS containing $2k+2=\frac{n+1}{2}$ vertices. Similarly for vertices $a_i,i\equiv 0,3(mod 4)$ consider the set $\{\{a_3,a_4\},\{a_7,a_8\},...,\{a_{4k+3},a_{1}\}\}$ which is a TDS containing $2k+2=\frac{n+1}{2}$ vertices. Hence $\mathpzc{d}_{td}(a_i)=\frac{n+1}{2}.$\\
    \textbf{Case 4:} $n\equiv2(mod 4).$\\
    Let $n=4k+2, k\in \mathbb{N}.$ For vertices of the form $a_i,i\equiv 1,2(mod 4)$ take the set $\{\{a_1,a_2\},\{a_5,a_6\},$\\$...,\{a_{4k+1},a_{4k+2}\}\}$ which is a TDS containing $2k+2=\frac{n}{2}+1$ vertices. Similarly for vertices $a_i,i\equiv 0,3(mod 4)$ consider the set $\{\{a_3,a_4\},\{a_7,a_8\},...,\{a_{4k-1},a_{4k}\},\{a_1,a_2\}\}$ which is a TDS containing $2k+2=\frac{n}{2}+1$ vertices. Hence $\mathpzc{d}_{td}(a_i)=\frac{n}{2}+1.$
\end{proof}
\begin{thm}\label{t6}
    A book graph $\mathbb{B}_n$ is a compliant graph and $$\mathpzc{d}_{td}(a_i)=\begin{cases}
     2 \quad &\text{ if  } a_i \text{ is center vertex}\\
    2n\quad &\text{ otherwise}   
   \end{cases}$$ 
\end{thm}
\begin{proof}
Let $\mathbb{B}_n= \mathbb{S}_{n+1}\times \mathbb{P}_2.$ Let the vertices of $\mathbb{S}_{n+1}$ be labeled as $b$ for the vertex with degree n, and all the remaining vertices as $b_1,b_2,...,b_n.$ Also, label the vertices of $\mathbb{P}_2$ as $a_1,a_2.$ Then in the book graph the center vertex is $(b,a_1)$ and $(b,a_2).$ Also $(b,a_1)$ is adjacent to $\{(b_1,a_1),(b_2,a_1),...,(b_n,a_1)\}.$ Similarly, $(b,a_2)$ is adjacent to $\{(b_1,a_2),(b_2,a_2),...,(b_n,a_2)\}.$ For the center vertices $(b,a_1)$ and $(b,a_2)$ the set $\{(b,a_1),(b,a_2)\}$ is a TDS containing the least number of vertices. Hence $\mathpzc{d}_{td}((b,a_1))=\mathpzc{d}_{td}((b,a_2))=2.$ For all other vertices $(b_i,a_j),i=1,2,...,n, j=1,2$ no MTDS contains $\{(b,a_1),(b,a_2)\}$ as subset. Hence set $\{\{(b_1,a_1),(b_1,a_2)\},\{(b_2,a_1),(b_2,a_2)\},...,\{(b_n,a_1),(b_n,a_2)\}\}$ with $2n$ vertices is a MTDS containing each vertex $(b_i,a_j)$ with least number of vertices. Hence $\mathpzc{d}_{td}((bi,a_j))=2n, i=1,2,...,n, j=1,2$ \end{proof}
\begin{prop}\label{pr1}
 For a Windmill graph $\mathbb{W}d(p,q)$, 
 $\mathpzc{d}_{td}(a)=2$
\end{prop}
\begin{proof}
 Consider the Windmill graph $\mathbb{W}d(p,q)$. Label the center vertex as $a$ and the vertices in the $i^{th}$ copy as $a_{i_1},a_{i_2},...,a_{i_{p-1}}.$ The center vertex dominates every other vertex. Hence for total domination, only one more vertex is required. Therefore, $\mathpzc{d}_{td}(a)=2.$ For any other vertex $a_{i_j},i=1,2,...,q\text{ and } j=1,2,...,p-1,$ the set $\{a_{i_j},a\}$ is a MTDS containing $a_{i_j}$ with least number of vertices. Hence  $\mathpzc{d}_{td}(a)=2$ for every vertex in $\mathbb{W}d(p,q).$ 
\end{proof}
\begin{rem}
    For a complete graph, complete bipartite graph, star graph, and wheel graph the TDD of every vertex is 2. Let the vertices of a wheel graph be labelled as $a,a_1,a_2,...,a_n,$ where $a$ is the center vertex. Then the set $\{a,a_i\}$ is a MTDS containing $a$ with the least cardinality. The same set can be taken as the MTDS containing $a_i$ with the least number of vertices. Hence the TDD of every vertex is 2.  
\end{rem}
\begin{figure}
    \centering
    \includegraphics[height=8cm,width=7.5cm]{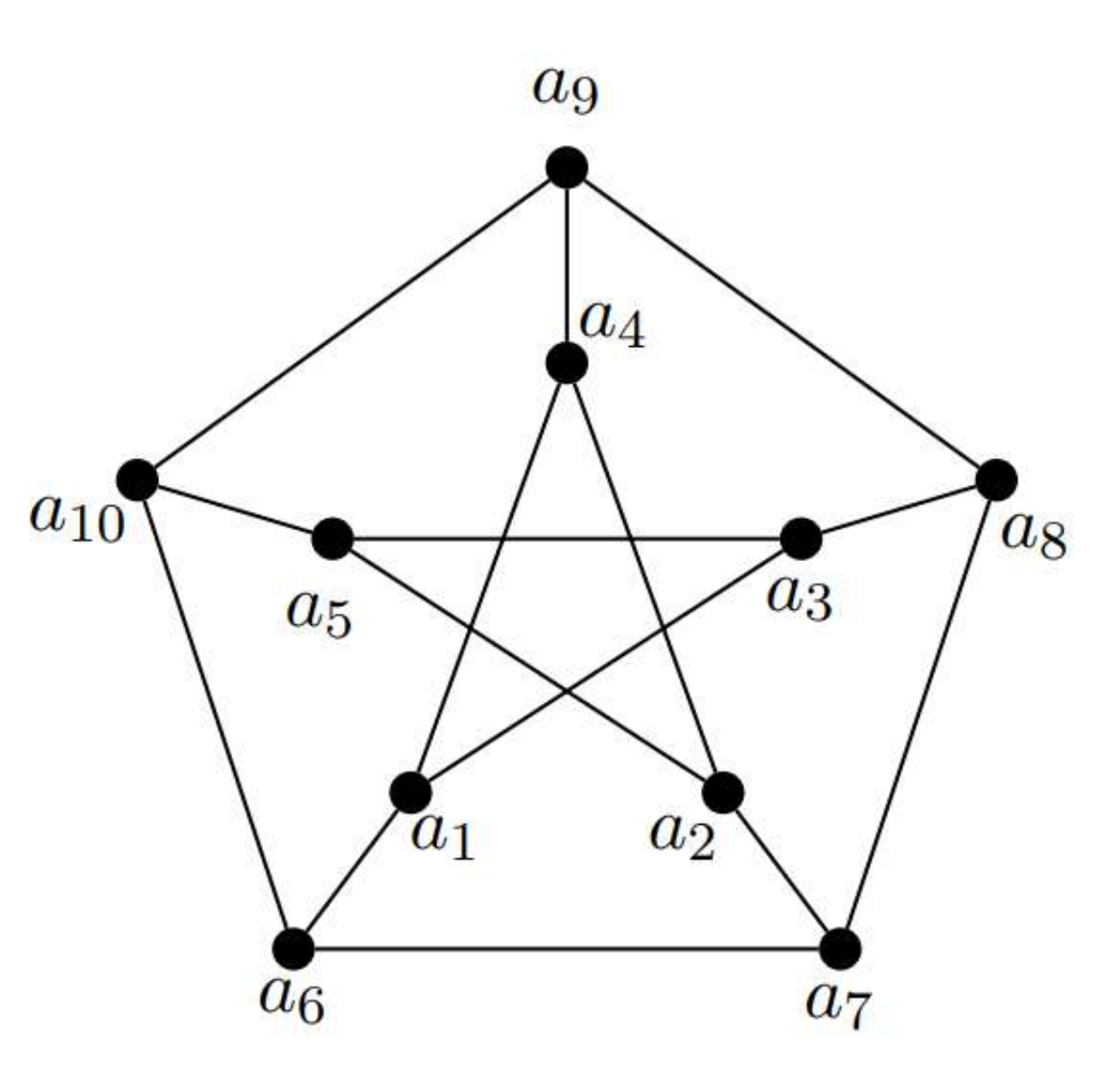}
    \caption{Petersen graph}
    \label{f3}
\end{figure}
\begin{exam}
    Consider the Petersen graph with labeling as in Figure \ref{f3}. The TDN of Petersen graph is 4. The sets $\{a_1,a_6,a_7,a_{10}\},\{a_2,a_6,a_7,a_8\},\{a_3,a_7,a_8,a_9\},\{a_4,a_8,a_9,a_{10}\},$ and $\{a_5,a_6,a_9,a_{10}\}$ can be considered as MTDS containing each of the vertices with cardinality 4. Hence TDD of each vertex is 4. Similarly, the TDD of every vertex is 4 for Gr\"{o}tzsch graph and Herschel graph. 
\end{exam}
\begin{defn}
    A total domination regular (TDR) graph is a graph in which every vertex has equal TDD. 
\end{defn}
\begin{exam}
    Petersen graph, complete graph, $\mathbb{P}_{4k+2},k=0,1,2...$ are examples of TDR graphs.
\end{exam}
\begin{prop}\label{pr2}
    For a compliant graph $\mathbb{G},$ $\gamma_t\leq \mathpzc{d}_{td}(a)\leq \Gamma_t,$ for every vertex $a\in V(\mathbb{G}).$
\end{prop}
\begin{proof}
    The inequality follows directly from the definition of TDN, TDD and upper TDN. 
\end{proof}
\begin{prop}\label{pr3}
    For a compliant graph $\mathbb{G},$ $\mathpzc{d}_{td}(a)\geq \mathpzc{d}_d(a).$
\end{prop}
\begin{proof}
    Since a TDS is always a DS the result follows.
\end{proof}
\begin{rem}
    Any lower bound for $\gamma_t$ is a lower bound for the TDD of a vertex. For example, $\frac{n}{\Delta(\mathbb{G})}$, $rad(\mathbb{G})$ and $\frac{diam(\mathbb{G})+1}{2}$ are lower bounds for $\gamma_t.$ Hence these lower bounds hold for TDD also.
\end{rem}
\begin{rem}
    Let $\mathbb{G}$ be a graph and $e$ be an edge in $\mathbb{G}.$ Then $\gamma_t(\mathbb{G}\setminus e)\geq \gamma_t(\mathbb{G}). $
\end{rem}
\begin{prop}\label{pr4}
    Let $\mathbb{G}$ be a compliant graph and $\mathbb{H}$ be a compliant subgraph of $\mathbb{G}.$ Then $_{\mathbb{G}}\mathpzc{d}_{td}(a)\leq {_{\mathbb{H}}\mathpzc{d}}_{td}(a)$ for every vertex $a\in V(\mathbb{H}).$  
\end{prop}
\begin{cor}
 Let $\mathbb{T}$ be a spanning tree of a connected compliant graph $\mathbb{G},$ such that $\mathbb{T}$ is also compliant. Then $_{\mathbb{G}}\mathpzc{d}_{td}(a)\leq{_{\mathbb{T}}\mathpzc{d}}_{td}(a)$ for every vertex $a\in V(\mathbb{G}).$ 
\end{cor}
\begin{thm}\label{t7}
    Let $\mathbb{G}=\bigcup\limits_{i=1}^m \mathbb{G}_i$ be the disjoint union of compliant graphs $\mathbb{G}_1, \mathbb{G}_2,...,\mathbb{G}_m$, then $$_{\mathbb{G}}\mathpzc{d}_{td}(a)={_{\mathbb{G}_i}\mathpzc{d}_{td}}(a)+\sum\limits_{\substack{{j=1}\\ j\neq i}}^m\gamma_{t}(\mathbb{G}_j).$$ 
\end{thm}
\begin{proof}
    Let $\mathbb{G}=\bigcup\limits_{i=1}^m \mathbb{G}_i$ be the disjoint union of compliant graphs $\mathbb{G}_1, \mathbb{G}_2,...,\mathbb{G}_m$.
    Let $a\in V(\mathbb{G}_k).$ A MTDS of $\mathbb{G}$ is disjoint union of MTDS of each component $\mathbb{G}_i.$ Any MTDS in $\mathbb{G}_j,j\neq k$ is independent of the vertex $a.$ Hence for components other than $\mathbb{G}_k$ a minimum TDS can be taken as MTDS. Now, for component $\mathbb{G}_k,$ a MTDS containing $a$ with least cardinality is required. Hence the cardinality of such MTDS will be ${_{\mathbb{G}_i}\mathpzc{d}_{td}}.$ Therefore the MTDS of $\mathbb{G}$ containing $a$ is $\mathbb{D}\cup \bigcup\limits_{\substack{{j=i}\\j\neq k}}^m\mathbb{D}_j$ where $\mathbb{D}_j$ is the minimum TDS of component $\mathbb{G}_j$ and $\mathbb{D}$ is the MTDS containing $a$ with least cardinality in $\mathbb{G}_k.$ Hence 
    $$_{\mathbb{G}}\mathpzc{d}_{td}(a)={_{\mathbb{G}_i}\mathpzc{d}_{td}}(a)+\sum\limits_{\substack{{j=1}\\ j\neq i}}^m\gamma_{t}(\mathbb{G}_j).$$
\end{proof}
\begin{prop}\label{pr5}
    Let $\mathbb{G}=\mathbb{G}_1 + \mathbb{G}_2$ be the join of two graphs $\mathbb{G}_1$ and $\mathbb{G}_2.$ Then, $_{\mathbb{G}}\mathpzc{d}_{td}(a)=2$ for every vertex in $\mathbb{G}.$
\end{prop}
\begin{proof}
    Let $\mathbb{G}$ be the join of $\mathbb{G}_1$ and $\mathbb{G}_2$. Since each vertex of $\mathbb{G}_1$ is adjacent to every vertex of $\mathbb{G}_2$ one vertex from each $\mathbb{G}_1$ and $\mathbb{G}_2$ is sufficient for the total domination. Hence for any vertex $a$ the set $\{a,b\}$ where $b $ is any vertex from the component not containing $a$ is a required MTDS with the least cardinality. Therefore, $_{\mathbb{G}}\mathpzc{d}_{td}(a)=2$ for every vertex in $\mathbb{G}.$
\end{proof}
\begin{thm}\label{t8}
    Let $\mathbb{G}$ be a compliant graph. Consider the composition $\mathbb{H}=\mathbb{G}\circ \mathbb{K}_n.$ Let the vertices of $\mathbb{G}$ be $a_1,a_2,...,a_m$ and the vertices of $\mathbb{K}_n$ be $b_1,b_2,...,b_n.$  Then, 
    $$_{\mathbb{H}}\mathpzc{d}_{td}(a_i,b_j)= {_{\mathbb{G}}\mathpzc{d}_{td}}(a_i).$$
\end{thm}
\begin{proof}
     Any vertex of $\mathbb{H}$ is of the form $(a_i,b_j), i=1,2,...,m$ and $j=1,2,...,n.$ By the definition of composition a vertex $(a_i,b_j) $ dominate all vertices of the form $(a_i,b_{j'}), j'=1,2,...,n$ and $(a_{i'},b_{j'}), j'=1,2,...,n,$ $a_{i'}$ are such that $a_i$ is adjacent to $a_{i'}$ in $\mathbb{G}.$ Hence any vertex $(a_i,b_j)$ is not adjacent to a vertex $(a,b)$ if $a_i$ is not adjacent to $a$ in $\mathbb{G}.$ Suppose that a MTDS containing $a_i$ in $\mathbb{G}$ is $\{a_i,c_1,c_2,...c_k\}, $ $c_1,c_2,...,c_k\in \{a_1,a_2,...,a_m\}\setminus \{a_i\}.$ Then $\mathbb{D}=$ $\{(a_i,b_j),(c_1,b_1),...,(c_k,b_1)\}$ is a TDS containing $(a_i,bj).$ Also $\mathbb{D}$ is minimal. Since if $\mathbb{D}\setminus \{(c_{k'},b_1)\}$ is a TDS, then $\{a_i,c_1,...c_k\}\setminus \{c_{k'}\}$ is also a TDS which is not possible. By the similar argument if $\mathbb{D}'$ is a MTDS containing $(a_i,b_j)$ having less number vertices than $\mathbb{D}$, then collecting all the first coordinate vertices of $\mathbb{D}',$ a MTDS of $\mathbb{G}$ is obtained containing $a_i$ with less number of vertices than that in $\{a_i,c_1,c_2,...c_k\}$ which is also not possible. Hence $\mathbb{D}=$ $\{(a_i,b_j),(c_1,b_1),...,(c_k,b_1)\}$ is a required set. Therefore $$_{\mathbb{H}}\mathpzc{d}_{td}(a_i,b_j)= _{\mathbb{G}}\mathpzc{d}_{td}(a_i).$$
\end{proof}
\begin{thm}\label{t9}
    Let $\mathbb{G}$ be a non-compliant graph. Then $\mathbb{H}=\mathbb{G}\circ \mathbb{K}_n$ is a compliant graph and $_{\mathbb{H}}\mathpzc{d}_{td}(a_i,b_j)\leq 2 _{\mathbb{G}}\mathpzc{d}_{d}(a_i).$
\end{thm}
\begin{proof}
    Let $\mathbb{G}$ be a non-compliant graph and let $\mathpzc{d}_d(a_i)$ be the domination degree of vertex $a_i$ in $\mathbb{G}.$ Suppose that the set $\mathbb{D}'=\{a_i,c_1,c_2,...,c_k\},$ $c_1,c_2,...,c_k\in \{a_1,a_2,...,a_m\}\setminus \{a_i\}$ be the MDS considered for the DD of $a_i.$ Then $\mathbb{D}=\{(a_i,b_1),(a_i,b_2),(c_1,b_1),(c_1,b_2),(c_2,b_1),(c_2,b_2),...,(c_k,b_1),(c_k,b_2)\}$ is a TDS with $2|\mathbb{D}'|$ vertices in $\mathbb{H}.$ Moreover, $\mathbb{D}$ is a MTDS since $\mathbb{D}'$ is MDS and if any vertex $(u_i,b_j)\in \mathbb{D}$ is deleted from $\mathbb{D}$ then $(u_i,b_k), k\neq j$ is not dominated. Hence each vertex is contained in at least one MTDS. Therefore $\mathbb{H}$ is compliant. Also, $_{\mathbb{H}}\mathpzc{d}_{td}(a_i,b_j)\leq 2 _{\mathbb{G}}\mathpzc{d}_{d}(a_i).$
\end{proof}
\begin{rem}
 In Theorem \ref{t9} non-compliant graphs are considered. In non-compliant graphs both compliant and non-compliant vertices exists. The bound for TDD of a vertex obtained in Theorem \ref{t9} is valid for all vertices, but for the compliant vertices in $\mathbb{G}$ the exact TDD is the same as in Theorem \ref{t8}.
\end{rem}
\begin{thm}\label{t10}
 Let $\mathbb{H}=\mathbb{G}_1\odot \mathbb{G}_2$ be the corona of a graph $\mathbb{G}_1$ and a compliant graph $\mathbb{G}_2.$ Then $\mathbb{H}$ is compliant. Let the vertices of $\mathbb{G}_1$ be $a_1,a_2,...,a_n$ and the vertices of $\mathbb{G}_2$ be $b_1,b_2,...,b_m $.  Let the vertices in the $i^{th}$ copy be $b_{i1},b_{i2},...,b_{im}, i=1,2,...,n.$ Then, 
    $$_{\mathbb{G}}\mathpzc{d}_{td}(a)=
     \begin{cases}
         n \, &\text{ if } a=a_i, i=1,2,...,n\\
         _{\mathbb{G}_2}\mathpzc{d}_{td}(a)+(n-1) \, &\text{ if } a=b_{ij}, i=1,2,...,n \text{ and } j=1,2,...,m.
     \end{cases}
    $$  
\end{thm}
\begin{proof}
     Let $\mathbb{H}=\mathbb{G}_1\odot \mathbb{G}_2$ be the corona of a graphs $\mathbb{G}_1$ and a compliant graph $\mathbb{G}_2.$ Any vertex of the form $a_i,i=1,2,...,n$ dominates all vertices in the $i^{th}$ copy and also the vertices which are adjacent to $a_i$ in $\mathbb{G}_1.$ To dominate any vertex of the form $b_{ij}$ either vertex $a_i$ or any other vertex of the form $b_{ij'}$ such that $b_j$ is adjacent to $b_{j'}$ in $\mathbb{G}_2$ is required. Hence at least n vertices are required for dominating all vertices. The set $\{a_1,a_2,...,a_n\}$ is sufficient for total domination of $\mathbb{H}.$ Hence for vertices of the form $a_i,$ $_{\mathbb{H}}\mathpzc{d}_{td}(a_i)=n.$ Now, for vertices of the form $b_{ij},$ the set $\{b_{ij},a_1,a_2,...,a_n\}$ is not a MTDS. Hence for $b_{ij},$ the MTDS containing $b_{j}$ in $\mathbb{G}_2$ is considered. And the remaining $a_k,k\neq i$ is chosen. Suppose $\mathbb{D}$ is the MTDS containing $b_{j}$ with the least cardinality in $\mathbb{G}_2$ then $\mathbb{D}\cup \{a_1,a_2,...,a_n\}\setminus a_i$ is the required MTDS for $b_{ij}$ with the least cardinality in $\mathbb{H}.$ Therefore $_{\mathbb{H}}\mathpzc{d}_{td}(b_{ij})={_{\mathbb{G}_2}\mathpzc{d}}_{td}(b_j)+(n-1).$
\end{proof}
\begin{thm}\label{t11}
 Let $\mathbb{G}_1$ be any graph and $\mathbb{G}_2$ be a non-compliant graph. Let $\mathbb{H}$ be the corona of $\mathbb{G}_1$ and $\mathbb{G}_2$. Then, $\mathbb{H}$ is compliant and $_{\mathbb{H}}\mathpzc{d}_{td}(b_{ij})\leq 2+(n-1)\gamma_t(\mathbb{G}_2),$ $i=1,2,...,n \text{ and } j=1,2,...,m.$  
\end{thm}
\begin{proof}
  Let $\mathbb{G}_1$ be any graph and $\mathbb{G}_2$ be a non-compliant graph. Let $\mathbb{H}$ be the corona of $\mathbb{G}_1$ and $\mathbb{G}_2$. For the vertices of the form $a_i,i=1,2,...,n$ the MTDS and TDD are the same as in Theorem \ref{t10}. For the non-compliant vertices of the form $b_{ij}, i=1,2,...,n$ $j=1,2,...,m,$ it is not possible to obtain a MTDS containing $b_j$ in $\mathbb{G}_2.$ Hence the set considered as in Theorem \ref{t10} is not possible. Now, take any minimum TDS say $\mathbb{D}$ in $\mathbb{G}_2$. Let $\mathbb{D}_i$ be the minimum TDS of $i^{th}$ copy of $\mathbb{G}_2$ which is the copy of $\mathbb{D}.$ For all copies of $\mathbb{G}_2$ other than $i^{th}$ copy take $\mathbb{D}_j$ to totally dominate all vertices in the copy and vertices $a_j, j\neq i.$ And for the $i^{th}$ copy consider $\{b_{ij},a_i\}. $ Hence the set $\{a_i,b_{ij}\}\cup \bigcup\limits_{\substack{{j'=1}\\ j\neq i}}^n \mathbb{D}_{j'} $ is a MTDS containing $b_{ij}$.  For all compliant $b_{ijs}$ the MTDS and TDD are same as in Theorem \ref{t10}. But the set $\{a_i,b_{ij}\}\cup \bigcup\limits_{\substack{{j'=1}\\ j\neq i}}^n \mathbb{D}_{j'} $ is also a MTDS containing $b_{ij}.$ Hence $\mathbb{H}$ is compliant and $_{\mathbb{H}}\mathpzc{d}_{td}(b_{ij})\leq 2+(n-1)\gamma_t(\mathbb{G}_2),$ $i=1,2,...,n \text{ and } j=1,2,...,m.$ 
\end{proof}
\begin{thm}\label{t12}\cite{sub}
For complete graph $\mathbb{K}_n,$ $$\gamma_t(S(\mathbb{K}_n)=\ceil[\Big]{\frac{3n}{2}}-1.$$ \end{thm}
\begin{thm}\label{t13}
 For a complete graph $\mathbb{K}_n,$ $S(\mathbb{K}_n)$ is compliant and $\mathpzc{d}_{td}(a)=\ceil[\Big]{\frac{3n}{2}}-1,$ for every vertex in $S(\mathbb{K}_n).$
\end{thm}
\begin{proof}
    From Theorem \ref{t12} the TDN of $S(\mathbb{K}_n)$ is $\ceil[\Big]{\frac{3n}{2}}-1.$ Hence it is enough to prove that there exists a MTDS containing each of the vertex in $S(\mathbb{K}_n).$ In $S(\mathbb{K}_n)$ let the new vertex between vertex $a_i$ and $a_j$ be labeled as $a_{ij}(a_{ji})$. \\
    \textbf{Case 1:} n is even.\\
    For any vertex of the form $a_{ij}$ consider the set $\{a_{ij},a_i,a_1,a_{12},a_2,a_{23},a_3,...a_{i-1},a_{(i-1)(i+1)},a_{i+1},$\\$...,a_{(n-1)n},a_n\}$ $\{\text{For $i=1$ choose set } \{a_{1j},a_1,a_2,a_{23},a_{3},...,a_{(n-1)},a_{(n-1)n},a_n\} \}$ which consists of $\ceil[\Big]{\frac{3n}{2}}-1$  vertices. The given set is a MTDS containing each vertices of the form $a_{ij}$ and $a_i.$ Hence TDD of every vertex is $\ceil[\Big]{\frac{3n}{2}}-1.$\\
    \textbf{Case 2:} n is odd.\\
    For any vertex of the form $a_{ij}$ consider the set $\{a_{ij},a_i,a_1,a_{12},a_2,a_{23},a_3,...a_{i-1},a_{(i-1)(i+1)},a_{i+1},$\\$...,a_{(n-2)},a_{(n-2)(n-1)}a_{(n-1)},a_n,a_{nj}\}$ $\{\text{For $i=1$ choose set } \{a_{1j},a_1,a_2,a_{23},a_{3},...,a_{(n-2)(n-1)}a_{(n-1)},a_n,a_{nj}\} \}$ which consists of $\ceil[\Big]{\frac{3n}{2}}-1$  vertices. The given set is a MTDS containing each vertices of the form $a_{ij}$ and $a_i.$ Hence TDD of every vertex is $\ceil[\Big]{\frac{3n}{2}}-1.$
\end{proof}
\begin{thm}
    Let $\mathbb{H}$ be the degree splitting graph of $\mathbb{P}_n.$ Then,
    $\mathpzc{d}_{td}(a)=\begin{cases}
        3 \text{ or } 4 \text{ for }\, \, n=4\\
        4 \text{ \qquad for }\, \, n>4
    \end{cases}$
\end{thm}
\begin{proof}
    Consider path $\mathbb{P}_n$ with vertices $a_1,a_2,...,a_n.$ Let the newly added vertices be labeled as $b_1$ and $b_2.$ Then $a_1$ and $a_n$ are joined to $b_1$ and all vertices with degree 2 is joined to $b_2.$ \\
    \textbf{Case 1:} $n=4.$\\
    Take $\mathbb{P}_4$ with vertices $a_1,a_2,a_3,$ and $a_4.$ Vertices $a_1$ and $a_4$ are made adjacent to $b_1.$ Similarly $a_2$ and $a_3$ are made adjacent to $b_2.$ Now, consider sets $\{a_1,a_2,b_1\}$ and $\{a_3,a_4,b_1\}.$ These two sets are MTDS containing each of $a_1,a_2,a_3,a_4,b_1$ with the least cardinality. Hence TDD of these vertices is 3. Similarly, for $b_2$ either $a_2$ or $a_3$ is required for total domination. Suppose $a_2$ is taken, then the vertices for total domination are $a_4$ and $b_1$. Hence at least 4 vertices are required for total domination and $\{b_2,a_2,a_4,b_1\}$ is a MTDS containing $b_2$ with the least cardinality. Hence TDD of $b_2$ is 4. \\
    \textbf{Case 2:} $n> 4.$\\
    Let $\mathbb{P}_n$ be path with $n>4$ vertices. Let $a_i$ be any vertex with degree 2. Then the set $\{a_1,b_1,b_2,a_i\}$ is a MTDS containing each vertex with the least cardinality. Hence TDD of every vertex is 4.\\
    Hence the result.
\end{proof}
\section{Total domination index of a graph}
\begin{defn}
    Let $\mathbb{G}$ be a compliant graph. The total domination index (TDI) of $\mathbb{G}$ is the sum of TDD of vertices of $\mathbb{G},$ i.e,
    $$ TDI(\mathbb{G})=\sum_{a\in V(\mathbb{G})} \mathpzc{d}_{td}(a).$$
 \end{defn}
 \begin{exam}
     For the graph in Figure \ref{f2}, there are two vertices with TDD 4, and 4 vertices with TDD 2. Hence TDI is $2\times 4+4 \times 2=16.$
 \end{exam}
 From Proposition \ref{pr2} and \ref{pr3} the bounds for TDI can be obtained. 
 \begin{prop}\label{pr6}
  For a compliant graph $\mathbb{G},$ $\gamma_t(\mathbb{G})\leq \frac{TDI}{n}\leq \Gamma_t(\mathbb{G}).$   
 \end{prop}
 \begin{prop}\label{pr7}
 For a compliant graph $\mathbb{G},$ $TDI(\mathbb{G})\geq DI(\mathbb{G}).$     
 \end{prop}
 \begin{prop}\label{pr8}
 Let $\mathbb{G}_1\cong \mathbb{G}_2,$ then $TDI(\mathbb{G}_1)=TDI(\mathbb{G}_2).$    
 \end{prop}
 \begin{proof}
     The proof is clear from the fact that isomorphism preserves adjacency and hence the TDD of each vertex. 
 \end{proof}
\begin{prop}\label{pr9}
 For cycle $\mathbb{C}_n,$
 $$TDI(\mathbb{C}_n)=\begin{cases}
     \frac{n^2}{2} \quad &\text{ if  } n\equiv 0 (mod4)\\
      \frac{n(n+1)}{2} \quad &\text{ if  } n\equiv 1,3 (mod 4)\\
      \frac{n^2}{2}+n \quad &\text{ if  } n\equiv 2 (mod  4)   
   \end{cases}$$ 
\end{prop}
\begin{prop}\label{pr10}
 For a book graph $\mathbb{B}_n,$ $TDI(\mathbb{B}_n)=4(1+n^2).$   
\end{prop}
\begin{proof}
    For a book graph $\mathbb{B}_n,$ there are two center vertices with TDD 2, and remaining $2n$ vertices have TDD $2n.$ Hence $TDI(\mathbb{B}_n)=2\times 2+2n\times 2n=4(1+n^2).$ 
\end{proof}
\begin{prop}\label{pr11}
 For a Windmill graph $\mathbb{W}(p,q),$ $TDI(\mathbb{W}(p,q))=2(q(p-1)+1).$
\end{prop}
\begin{prop}\label{pr12}
For a complete graph $\mathbb{K}_n,$ $TDI(\mathbb{K}_n)=2n.$    
\end{prop}
\begin{prop}\label{pr13}
For a complete bipartite graph $\mathbb{K}_{m,n}$ $TDI(\mathbb{K}_n)=2(m+n).$    
\end{prop}
\begin{prop}\label{pr14}
For a star graph $\mathbb{K}_{1,n}$ and wheel graph $\mathbb{W}_n,$  $TDI(\mathbb{G})=2(1+n).$  
\end{prop}
\begin{thm}\label{t14}
    For path graph $\mathbb{P}_n,n\neq4,7$ and for $k=0,1,2,3,...$
    $$TDI(\mathbb{P}_n)=\begin{cases}
4(2k+1)(k+1)  &\text{ if } n=4k+2,k\geq0\\
8k^2+7k+2 & \text{ if } n=4k+1,k>0\\
8k^2+15k+6 & \text{ if } n=4k+3,k=0,k>1\\
2k(4k+1) & \text{ if } n=4k,k>1
    \end{cases}$$
\end{thm}
\begin{proof}
    Consider path $\mathbb{P}_n,n\neq 4,7.$\\
    \textbf{Case 1:} $n=4k+2, k\geq 0.$\\
   From Theorem \ref{t4}, for $n\equiv 2(mod 4)$ every vertex has TDD $\frac{n}{2}+1.$ Hence $TDI(\mathbb{P}_n)=n\times \frac{n}{2}+1= (4k+2)\times (\frac{4k+2}{2}+1)=4(2k+1)(k+1).$\\
   \textbf{ Case 2:} $n=4k+1, k> 0.$\\
   From Theorem \ref{t4}, for $n\equiv 1(mod 4)$ the TDD is $\frac{n+1}{2}+1$ for vertices $a_i$ such that $i\equiv1(mod 4)$ and for $i\not\equiv1(mod4) $ TDD is $\frac{n+1}{2}.$ For $n=4k+1,$ there are $k+1$ vertices with $i\equiv 1(mod4)$ and $3k$ vertices with $i\not\equiv1(mod4).$ Hence $TDI(\mathbb{P}_n)=(k+1)\times (\frac{4k+1+1}{2}+1)+3k(\frac{4k+1+1}{2})=8k^2+7k+2.$\\
   \textbf{ Case 3:} $n=4k+3, k=0,k> 1.$\\
   From Theorem \ref{t4}, for $n\equiv 3(mod 4)$ the TDD is $\frac{n+1}{2}+1$ for vertices $a_i$ such that $i\equiv0(mod 4)$ and for $i\not\equiv1(mod4) $ TDD is $\frac{n+1}{2}.$ For $n=4k+3,$ there are $k$ vertices with $i\equiv 0(mod4)$ and $3k+3$ vertices with $i\not\equiv0(mod4).$ Hence $TDI(\mathbb{P}_n)=(k)\times (\frac{4k+3+1}{2}+1)+(3k+3)(\frac{4k+3+1}{2})=8k^2+15k+6.$\\
   \textbf{ Case 4:} $n=4k, k> 1.$\\
   From Theorem \ref{t4}, for $n\equiv 0(mod 4)$ the TDD is $\frac{n}{2}$ for vertices $a_i$ such that $i\equiv2,3(mod 4)$ and for $i\equiv0,1(mod4) $ TDD is $\frac{n}{2}+1.$ For $n=4k,$ there are $k+k$ vertices with $i\equiv 0,1(mod4)$ and $2k$ vertices with $i\equiv2,3(mod4).$ Hence $TDI(\mathbb{P}_n)=(2k)\times (\frac{4k}{2}+1)+2k(\frac{4k}{2})=2k(4k+1).$
\end{proof}
\begin{prop}
    Let $\mathbb{H}$ be a spanning compliant subgraph of $\mathbb{G}$, then $TDI(\mathbb{G})\leq TDI(\mathbb{H}).$
\end{prop}
\begin{cor}
   Let $\mathbb{T}$ be a compliant spanning tree of $\mathbb{G}$, then $TDI(\mathbb{G})\leq TDI(\mathbb{T}).$ 
\end{cor}
\section{Conclusion}
In conclusion, this paper introduces the concepts of compliant vertices and graphs, focusing on the total domination degree (TDD) and total domination index (TDI) across various graph classes. TDD is calculated for several types of graphs, with inequalities and operations like union, join, and composition explored. The study also examines subdivision and degree splitting. The TDI is introduced, and its values and bounds are determined for different graph types. This work lays the foundation for further research on domination parameters and their applications in graph theory. 

\vspace{2cm}
\textbf{\large{Acknowledgement}}\\
The first author gratefully acknowledges the financial support of the Council of Science and Industrial Research (CSIR), Government of India.\\
The authors would like to thank the DST, Government of India, for providing support to carry out this work under the scheme 'FIST' (No. SR/FST/MS-I/2019/40).

\end{document}